\documentclass{article}

\begin{document}
\title{ On the absolute matrix summability factors }
\author{H. S. \"OZARSLAN and T. ARI\\
Department of Mathematics, Erciyes University, 38039 Kayseri,
Turkey\\ E-mail:seyhan@erciyes.edu.tr and tkandefer@erciyes.edu.tr
\date{}} \maketitle
\begin{abstract}
In this paper, we have obtained a necessary and sufficient condition on $(\lambda_{n})$ for the series $\sum \lambda_{n} a_{n}$ to be $\left|A\right|_{k}$ summable, $k\geq1$, whenever $\sum a_{n}$ is $\left|A\right|$ summable. As a consequence we extend some known results of Sar{\i}g\"ol [2].
\end{abstract}
\bigskip
\noindent {\bf 1. Introduction} \\
\\Let $\sum a_{n}$ be a given infinite series with the partial sums
$\left( s_{n}\right) $, and let $A=(a_{nv})$ be a normal matrix,
i.e., a lower triangular matrix of nonzero diagonal entries. Then
$A$ defines the sequence-to-sequence transformation, mapping the
sequence $s=(s_{n})$ to $As=\left(A_{n}(s)\right)$, where
\begin{eqnarray}
A_{n}(s)=\sum_{v=0}^{n}a_{nv}s_{v},\quad n=0,1,...
\end{eqnarray}
The series $\sum a_{n}$ is said to be summable $\left| A\right|
_{k}\,,k\geq 1$, if (see [3])
\begin{eqnarray}
\sum_{n=1}^{\infty}n^{k-1}\
\left|\bar{\Delta}A_{n}(s)\right|^{k}<\infty ,
\end{eqnarray}
where
\begin{eqnarray*}
\bar{\Delta}A_{n}(s)=A_{n}(s)-A_{n-1}(s)
\end{eqnarray*}
and it is said to be
$\left|R,p_{n}\right|_{k}$ summable (see [5])if (2) holds when $A$ is a
Riesz matrix.
\bigskip
\\By a Riesz matrix we mean one such that
\begin{eqnarray*}
a_{nv}=\frac{p_{v}}{P_{n}}, \quad for \quad 0\leq v\leq n,\quad
and \quad a_{nv}=0 \quad for \quad v>n,
\end{eqnarray*}
where $(p_{n})$ is a sequence of positive real numbers such that
\begin{eqnarray*}
P_{n}=\sum_{v=0}^{n}p_{v}\rightarrow\infty,\quad
(n\rightarrow\infty), \quad \left(P_{-i}=p_{-i}=0,\quad i\geq 1
\right).
\end{eqnarray*}
Sar{\i}g\"ol [2] has proved the following theorem for
$\left|R,p_{n}\right|_{k}$ summability method.
\bigskip
\\ \noindent {\bf Theorem A.}
Suppose that $(p_{n})$ and $(q_{n})$ are positive sequences with $P_{n}\rightarrow\infty$ and $Q_{n}\rightarrow\infty$ as $n\rightarrow\infty$. Then $\sum  a_{n} \lambda_{n}$ is summable $\left| R,q_{n}\right| _{k}$, $k\geq1$, whenever $\sum a_{n}$ is summable $\left| R,p_{n}\right|$, if and only if
\begin{eqnarray*}
\textbf{(a)} \ \ \lambda_{n}=O\left\{ n^{\frac{1}{k}-1} \frac{q_{n}P_{n}}{p_{n}Q_{n}}\right\},
\end{eqnarray*}
\begin{eqnarray}
\textbf{(b)} \ \ W_{n}\triangle \left( Q_{n-1} \lambda_{n} \right) =O \left( \frac{p_{n}}{P_{n}} \right),
\end{eqnarray}
\begin{eqnarray*}
\textbf{(c)} \ \ Q_{n}\lambda_{n+1}W_{n}=O(1),
\end{eqnarray*}
where, provided that
\begin{eqnarray*}
W_{n}=\left\{ \sum_{v=n+1}^{\infty}v^{k-1} \left( \frac{q_{v}}{Q_{v}Q_{v-1}} \right)^{k} \right\}^{\frac{1}{k}}<\infty.
\end{eqnarray*}
\bigskip
\\ \noindent {\bf Lemma. ([4])}
$A=(a_{nv})\in (l_{1},l_{k})$ if and only if
\begin{eqnarray}
\sup_{v} \sum_{n=1}^{\infty} |a_{nv}|^{k}<\infty
\end{eqnarray}
for the cases $1\leq k<\infty$, where $(l_{1},l_{k})$ denotes the set of all matrices $A$ which map $l_{1}$ into $l_{k}=\{ x=(x_{n})\ :\ \sum |x_{n}|^{k} <\infty\}$.

\bigskip

\noindent {\bf 2. The main result.} The aim of this paper is to
generalize Theorem $A$ for absolute matrix summability. Before
stating the main theorem we must first introduce some further
notations.
\\Given a normal matrix $A=(a_{nv})$, we associate two lover
semimatrices $\bar{A}=(\bar{a}_{nv})$ and $\hat{A}=(\hat{a}_{nv})$
as follows:
\begin{eqnarray}
\bar{a}_{nv}=\sum_{i=v}^{n}a_{ni},\quad n,v=0,1,...
\end{eqnarray}
and
\begin{eqnarray}
\hat{a}_{00}=\bar{a}_{00}=a_{00},\quad
\hat{a}_{nv}=\bar{a}_{nv}-\bar{a}_{n-1,v}\quad n=1,2,...
\end{eqnarray}
It may be noted that $\bar{A}$ and $\hat{A}$ are the well-known
matrices of series-to-sequence and series-to-series
transformations, respectively. Then, we have
\begin{eqnarray}
A_{n}(s) & = &
\sum_{v=0}^{n}a_{nv}s_{v}=\sum_{v=0}^{n}a_{nv}\sum_{i=0}^{v}a_{i} \nonumber \\
 & = & \sum_{i=0}^{n}a_{i}\sum_{v=i}^{n}a_{nv}= \sum_{i=0}^{n}\bar{a}_{ni}a_{i}
\end{eqnarray}
and
\begin{eqnarray}
\bar{\Delta}A_{n}(s)& = & \sum_{i=0}^{n}\bar{a}_{ni}a_{i}-
\sum_{i=0}^{n-1}\bar{a}_{n-1,i}a_{i} \nonumber \\ & = &
\bar{a}_{nn} a_{n}+\sum_{i=0}^{n-1}(\bar{a}_{ni}-\bar{a}_{n-1,i})a_{i}
\nonumber \\ & = &
\hat{a}_{nn}a_{n}+\sum_{i=0}^{n-1}\hat{a}_{ni}a_{i}=\sum_{i=0}^{n}\hat{a}_{ni}a_{i}.
\end{eqnarray}
If $A$ is a normal matrix, then $A'=(a'_{nv})$ will denote the
inverse of $A$. Clearly, if $A$ is normal then
$\hat{A}=(\hat{a}_{nv})$ is normal and it has two-sided inverse
$\hat{A}'=(\hat{a}'_{nv})$, which is also normal (see [1]).
\\Now we shall prove the following theorem.

\bigskip

\noindent {\bf Theorem.} Let $k\geq1$, $A=(a_{nv})$ and
$B=(b_{nv})$ be two positive normal matrices. In order that $\sum  a_{n} \lambda_{n}$ is summable $\left| B \right| _{k}$, whenever $\sum a_{n}$ is summable $\left| A \right|$ it is necessary that
\begin{eqnarray}
|\lambda_{n}|=O \left\{ n^{\frac{1}{k}-1} \frac{a_{nn}}{b_{nn}} \right\},
\end{eqnarray}
\begin{eqnarray}
\sum_{n=v+1}^{\infty}n^{k-1}|\Delta_{v}(\hat{b}_{nv }\lambda_{v})|^{k}=O(a_{vv})^{k},
\end{eqnarray}
\begin{eqnarray}
\sum_{n=v+1}^{\infty}n^{k-1} |\hat{b}_{n,v+1} \lambda_{v+1}|^{k}=O(1),
\end{eqnarray}
\begin{eqnarray}
a_{n-1,v} \geq a_{nv},\quad for \quad n\geq v+1,
\end{eqnarray}
\begin{eqnarray}
\bar{a}_{n0}=1, \quad  n=0,1,2,... \ .
\end{eqnarray}
Then (9)-(11) and
\begin{eqnarray}
\bar{b}_{n0}=1, \quad n=0,1,2,...,
\end{eqnarray}
\begin{eqnarray}
a_{nn}-a_{n+1,n}=O(a_{nn} \ a_{n+1,n+1}),
\end{eqnarray}
\begin{eqnarray}
\sum_{v=r+2}^{n}\left |\hat{b}_{nv}\right |\left
|\hat{a}'_{vr}\lambda_{v}\right |=O(  \frac{b_{nn}}{a_{nn}}| \lambda_{n}  |)
\end{eqnarray}
are also sufficient.\\
\\It should be noted that if we take $a_{nv}=\frac{p_{v}}{P_{n}}$
and $b_{nv}=\frac{q_{v}}{Q_{n}}$, then we get Theorem A.\\
\\ \noindent {\bf Proof of the theorem.}\\
\noindent {\bf  Necessity.} Let $(x_n)$ and $(y_n)$ denote $A$-transform and $B$-transform of
the series $\sum a_n$ and $\sum a_n\lambda_{n}$, respectively. Then, by (7) and (8), we have
\begin{eqnarray}
\overline \Delta x_n=\sum_{v=0}^{n}\hat{a}_{nv} a_v
 \ and \   \overline \Delta y_n=\sum_{v=0}^{n}\hat{b}_{nv} a_v \lambda_{v}.
\end{eqnarray}
For $k \geq 1$, we define
\begin{eqnarray*}
A=\left\{( a_{i}) : \sum a_{i}\ is \ summable \ |A| \right\},
\end{eqnarray*}
\begin{eqnarray*}
B=\left\{ (a_{i} \lambda_{i}) : \sum a_{i}\lambda_{i} \ is\ summable \ |B|_{k}\right\}.
\end{eqnarray*}
Then it is routine to verify that these are BK-spaces, if normed
by
\begin{eqnarray}
\left\| X\right\| =\left\{\sum_{n=0}^{\infty
} \mid{ \overline \Delta
x_{n}}\mid \right\}
\end{eqnarray}
and
\begin{eqnarray}
\left\| Y\right\| =\left\{\sum_{n=0}^{\infty
} n^{k-1} \mid{  \overline \Delta
y_{n}}\mid ^{k}\right\} ^{\frac{1}{k}}
\end{eqnarray}
respectively.\\
Since $\sum a_{n}$ is summable $|A|$ implies
$\sum a_{n}\lambda_{n}$ is summable $|B|_{k}$, by the
hypothesis of the theorem,
\begin{eqnarray*}
\left\| X\right\| <\infty \Rightarrow \left\| Y\right\| <\infty .
\end{eqnarray*}
Now consider the inclusion map c: A$\rightarrow $B defined by
c(x)=x. This is continous, which is immediate as A and B are
BK-spaces. Thus there exists a constant \ M such that
\begin{eqnarray}
\left\| Y\right\| \leq M\,\left\| X\right\|.
\end{eqnarray}
By applying (17) to $a_{v}=e_{v}-e_{v+1}$ ( $e_{v}$ is the v-th
coordinate vector), we have
\[\overline{\Delta} x_{n}=\left\{\begin{array}{cl}
                                   0                          &, \mbox{    if  $n<v$}\\
                                   \hat{a}_{nv}           &, \mbox{    if $n=v$}\\
                                   \Delta_{v}\hat{a}_{nv} &, \mbox{    if $n>v$}
                                   \end{array}
        \right. \]
and
\[\overline{\Delta} y_{n}=\left\{\begin{array}{cl}
                                   0                          &, \mbox{    if $n<v$}\\
                                   \hat{b}_{nv} \lambda_{v}           &, \mbox{    if $n=v$}\\
                                  \Delta_{v}( \hat{b}_{nv} \lambda_{v}) &, \mbox{    if $n>v$}.
                                   \end{array}
        \right. \]
So (18) and (19) give us
\begin{eqnarray*}
\left\| X\right\| =\left\{ a_{vv}+ \sum_{n=v+1}^{\infty} \mid{\Delta_{v}\hat{a}_{nv}
}\mid \right\}
\end{eqnarray*}
and
\begin{eqnarray*}
\left\| Y\right\| =\left\{v
^{k-1}b_{vv}\mid{\lambda_{v}}\mid^{k}+ \sum_{n=v+1}^{\infty} n ^{k-1}\mid{\Delta_{v}\left(\hat{b}_{nv}\lambda_{v}\right)
}\mid^{k}\right\} ^{\frac{1}{k}}.
\end{eqnarray*}
Hence it follows from (20) that\\
\begin{eqnarray*}
v^{k-1}b_{vv}\mid{\lambda_{v}}\mid^{k}+\sum_{n=v+1}^{\infty}n^{k-1}\mid{\Delta_{v}\hat{b}_{nv}\lambda_{v}}\mid^{k}
& \leq & M^{k}a_{vv}^{k}+ M^{k}\sum_{n=v+1}^{\infty}
\mid{\Delta_{v}\hat{a}_{nv}}\mid^{k}.
\end{eqnarray*}
Using (12), we can find
\begin{eqnarray*}
v^{k-1}b_{vv}\mid{\lambda_{v}}\mid^{k}+\sum_{n=v+1}^{\infty}n^{k-1}\mid{\Delta_{v}(\hat{b}_{nv}\lambda_{v})}\mid^{k} =
O\left\{ a_{vv}^{k}\right\}.
\end{eqnarray*}
The above inequality will be true iff each term on the left hand
side is  $O\left\{ a_{vv}^{k}\right\}$. Taking the first term,\\
\begin{eqnarray*}
v^{k-1}b_{vv}\mid{\lambda_{v}}\mid^{k}=
O\left\{ a_{vv}^{k}\right\}
\end{eqnarray*}
then
\begin{eqnarray*}
\mid{\lambda_{v}}\mid=
O\left\{v^{\frac{1}{k}-1} \frac{a_{vv}}{b_{vv}} \right\}
\end{eqnarray*}
which verifies that (9) is necessary.\\
Using the second term we have,
\begin{eqnarray*}
\sum_{n=v+1}^{\infty}n^{k-1}\mid{\Delta_{v}(\hat{b}_{nv}\lambda_{v})}\mid^{k}=
O\left\{ \mid{a_{vv}}\mid^{k}\right\}
\end{eqnarray*}
which is condition (10).\\
Now if we apply (17) to $a_{v}=e_{v+1}$, we have,
\[\overline{\Delta} x_{n}=\left\{\begin{array}{cl}
                                   0                   &, \mbox{    if $n\leq v$}\\
                                   \hat{a}_{n,v+1} &, \mbox{    if $n>v$}
                                   \end{array}
        \right. \]
and
\[\overline{\Delta} y_{n}=\left\{\begin{array}{cl}
                                   0                   &, \mbox{    if $n\leq v$}\\
                                   \hat{b}_{n,v+1}\lambda_{v+1} &, \mbox{    if $n>v$}
                                   \end{array}
        \right. \]
respectively.\\
Hence
\begin{eqnarray*}
\left\| X\right\|
=\left\{ \sum_{n=v+1}^{\infty}\mid{\hat{a}_{n,v+1}}\mid\right\}
,
\end{eqnarray*}
\begin{eqnarray*}
\left\| Y\right\| =\left\{
\sum_{n=v+1}^{\infty}n^{k-1}\mid{\hat{b}_{n,v+1}\lambda_{v+1}}\mid^{k}\right\}
^{\frac{1}{k}}.
\end{eqnarray*}
Hence it follows from (20) that
\begin{eqnarray*}
\sum_{n=v+1}^{\infty}n^{k-1}\mid{\hat{b}_{n,v+1}\lambda_{v+1}}\mid^{k} \leq M^{k} \left\{ \sum_{n=v+1}^{\infty
}\mid{\hat{a}_{n,v+1}}\mid\right\}^{k}.
\end{eqnarray*}
Using (13) we can find
\begin{eqnarray*}
\sum_{n=v+1}^{\infty}n^{k-1}\mid{\hat{b}_{n,v+1}\lambda_{v+1}}\mid^{k}=O(1)
\end{eqnarray*}
which is condition (11).\\

\bigskip

\noindent {\bf Sufficiency.} We use the notations of necessity. Then
\begin{eqnarray}
\overline \Delta x_n=\sum_{v=0}^{n}\hat{a}_{nv} a_v
\end{eqnarray}
which implies
\begin{eqnarray}
a_{v}=\sum_{r=0}^{v}\hat{a}'_{vr}\ \overline{\Delta}x_{r}.
\end{eqnarray}
In this case
\begin{eqnarray*}
\bar{\Delta}y_{n}=\sum_{v=0}^{n}\hat{b}_{nv}a_{v}\lambda_{v}=\sum_{v=0}^{n}\hat{b}_{nv}\lambda_{v}\
\sum_{r=0}^{v}\hat{a}'_{vr}\bar{\Delta}x_{r}.
\end{eqnarray*}
On the other hand, since
\begin{eqnarray*}
\hat{b}_{n0}=\bar{b}_{n0}-\bar{b}_{n-1,0}
\end{eqnarray*}
by (14), we have
\begin{eqnarray}
\bar{\Delta}y_{n} & = & \sum_{v=1}^{n}\hat{b}_{nv}\lambda_{v}\{\sum_{r=0}^{v}\hat{a}'_{vr}\ \bar{\Delta}x_{r}\}\nonumber \\
                         & = & \sum_{v=1}^{n}\hat{b}_{nv}\lambda_{v}\{\hat{a}'_{vv}\ \bar{\Delta}x_{v}+\hat{a}'_{v,v-1}\ \bar{\Delta}x_{v-1}+\sum_{r=0}^{v-2}\hat{a}'_{vr}\ \bar{\Delta}x_{r}\}\nonumber \\
                         & = & \sum_{v=1}^{n}\hat{b}_{nv}\lambda_{v}\ \hat{a}'_{vv}\ \bar{\Delta}x_{v}+\sum_{v=1}^{n}\hat{b}_{nv}\lambda_{v}\ \hat{a}'_{v,v-1}\ \bar{\Delta}x_{v-1}+\sum_{v=1}^{n}\hat{b}_{nv}\lambda_{v}\sum_{r=0}^{v-2}\hat{a}'_{vr}\ \bar{\Delta}x_{r}\nonumber \\
                         & = & \hat{b}_{nn}\lambda_{n}\ \hat{a}'_{nn}\ \bar{\Delta}x_{n}+\sum_{v=1}^{n-1}(\hat{b}_{nv}\lambda_{v}\ \hat{a}'_{vv}+\ \hat{b}_{n,v+1}\lambda_{v+1}\ \hat{a}'_{v+1,v})\ \bar{\Delta}x_{v} \nonumber \\
                         & \quad & +\sum_{r=0}^{n-2}\bar{\Delta}x_{r}\sum_{v=r+2}^{n}\hat{b}_{nv}\lambda_{v}\ \hat{a}'_{vr}.
\end{eqnarray}
By considering the equality
\begin{eqnarray*}
\sum_{k=v}^{n}\hat{a}'_{nk}\hat{a}_{kv}=\delta_{nv}
\end{eqnarray*}
where $\delta_{nv}$ is the Kronocker delta, we have that
\begin{eqnarray*}
\hat{b}_{nv}\lambda_{v} \ \hat{a}'_{vv}+\hat{b}_{n,v+1}\lambda_{v+1} \ \hat{a}'_{v+1,v} & =
& \frac{\hat{b}_{nv}\lambda_{v}}{\hat{a}_{vv}}+\hat{b}_{n,v+1}\lambda_{v+1}\
(-\frac{\hat{a}_{v+1,v}}{\hat{a}_{vv}\
\hat{a}_{v+1,v+1}})\nonumber \\    & = &
\frac{\hat{b}_{nv}\lambda_{v}}{a_{vv}}-\frac{\hat{b}_{n,v+1}\lambda_{v+1}\
(\bar{a}_{v+1,v}-\bar{a}_{v,v})}{a_{vv}\ a_{v+1,v+1}}\nonumber \\
& = & \frac{\hat{b}_{nv}\lambda_{v}}{a_{vv}}-\frac{\hat{b}_{n,v+1}\lambda_{v+1}\
(a_{v+1,v+1}+a_{v+1,v}-a_{vv})}{a_{vv}\ a_{v+1,v+1}}\nonumber \\ &
= & \frac{\Delta_{v}\left(\hat{b}_{nv}\lambda_{v}\right)}{a_{vv}}+\hat{b}_{n,v+1}\lambda_{v+1}\
\frac{a_{vv}-a_{v+1,v}}{a_{vv}\ a_{v+1,v+1}}
\end{eqnarray*}
and so
\begin{eqnarray*}
\bar{\Delta}y_{n} & = & \frac{b_{nn}\lambda_{n}}{a_{nn}}\
\bar{\Delta}x_{n}+\sum_{v=1}^{n-1}\
\frac{\Delta_{v}\left(\hat{b}_{nv}\lambda_{v}\right)}{a_{vv}}\
\bar{\Delta}x_{v}+\sum_{v=1}^{n-1}\hat{b}_{n,v+1}\lambda_{v+1}\
\frac{a_{vv}-a_{v+1,v}}{a_{vv}\ a_{v+1,v+1}}\
\bar{\Delta}x_{v} \nonumber \\    & + & \sum_{r=0}^{n-2}\bar{\Delta}x_{r}\sum_{v=r+2}^{n}\hat{b}_{nv}\lambda_{v}\
\hat{a}'_{vr}.
\end{eqnarray*}
Let
\begin{eqnarray*}
T_{n}(1)=\frac{b_{nn}\lambda_{n}}{a_{nn}}\
\bar{\Delta}x_{n}+\sum_{v=1}^{n-1}\
\frac{\Delta_{v}\left(\hat{b}_{nv}\lambda_{v}\right)}{a_{vv}}\
\bar{\Delta}x_{v}+\sum_{v=1}^{n-1}\hat{b}_{n,v+1}\lambda_{v+1}\
\frac{a_{vv}-a_{v+1,v}}{a_{vv}\ a_{v+1,v+1}}\
\bar{\Delta}x_{v},
\end{eqnarray*}
\begin{eqnarray*}
T_{n}(2)=\sum_{r=0}^{n-2}\bar{\Delta}x_{r}\sum_{v=r+2}^{n}\hat{b}_{nv}\lambda_{v}\
\hat{a}'_{vr}.
\end{eqnarray*}
Since
\begin{eqnarray*}
\left|T_{n}(1)+T_{n}(2)\right| ^{k} \leq 2^{k}
\left(\left|T_{n}(1)\right|^{k}+\left|T_{n}(2)\right|^{k}\right)
\end{eqnarray*}
to complete the proof of theorem, it is sufficient to show that
\begin{eqnarray*}
\sum_{n=1}^{\infty}n^{k-1}\left| T_{n}(i) \right| ^{k}<\infty
\quad for \quad i=1,2.
\end{eqnarray*}
Then
\begin{eqnarray*}
\overline{T_{n}(1)} & = & n^{1-\frac{1}{k}} \ T_{n}(1)\nonumber \\    & = & n^{1-\frac{1}{k}}\frac{b_{nn}\lambda_{n}}{a_{nn}}\
\bar{\Delta}x_{n}+n^{1-\frac{1}{k}} \sum_{v=1}^{n-1}\
\frac{\Delta_{v}\left(\hat{b}_{nv}\lambda_{v}\right)}{a_{vv}}\
\bar{\Delta}x_{v}+n^{1-\frac{1}{k}}\sum_{v=1}^{n-1}\hat{b}_{n,v+1}\lambda_{v+1}\
\frac{a_{vv}-a_{v+1,v}}{a_{vv}\ a_{v+1,v+1}}\
\bar{\Delta}x_{v}\nonumber \\    & = & \sum_{v=1}^{\infty} c_{nv} \bar{\Delta}x_{v}
\end{eqnarray*}
where
\[c_{nv}=\left\{\begin{array}{cl}
                                   n^{1-\frac{1}{k}}\left( \frac{\Delta_{v}\left({b}_{nv}\lambda_{v}\right)}{a_{vv}}+ \hat{b}_{n,v+1}\lambda_{v+1}\
\frac{a_{vv}-a_{v+1,v}}{a_{vv}\ a_{v+1,v+1}}   \right)  &, \mbox{    if $1\leq v\leq n-1$}\\
                                   n^{1-\frac{1}{k}} \frac{{b}_{nn}\lambda_{n}}{a_{nn}} &, \mbox{    if $v=n$}\\
                                   0 &, \mbox{    if $v>n.$}
                                   \end{array}
        \right. \]
Now
\begin{eqnarray*}
\sum |\overline{T_{n}(1)}|^{k}<\infty \ \ \textmd{whenever} \ \ \sum |\bar{\Delta}x_{n}|<\infty
\end{eqnarray*}
is equivalently
\begin{eqnarray}
\sup_{v} \sum _{n=1}^{\infty} |c_{nv}|^{k}<\infty
\end{eqnarray}
by Lemma. But (24) is equivalent to
\begin{eqnarray}
\sum _{n=v}^{\infty} |c_{nv}|^{k} & = & O(1)\left\{n^{1-\frac{1}{k}}| \frac{{b}_{nn}\lambda_{n}}{a_{nn}}|^{k}+ \sum_{n=v+1}^{\infty}  n^{1-\frac{1}{k}}\left| \frac{\Delta_{v}\left(\hat{b}_{nv}\lambda_{v}\right)}{a_{vv}}+ \hat{b}_{n,v+1}\lambda_{v+1}\
\frac{a_{vv}-a_{v+1,v}}{a_{vv}\ a_{v+1,v+1}}   \right|^{k}\right\} \nonumber \\    & = & O(1) \ \ as \ \ v \rightarrow \infty.
\end{eqnarray}
Finally
\begin{eqnarray*}
\sum_{n=2}^{\infty}n^{k-1}\left| T_{n}(2) \right|^{k} & = &
\sum_{n=2}^{\infty}n^{k-1} \left| \sum_{r=0}^{n-2}
\bar{\Delta}x_{r}\sum_{v=r+2}^{n}\hat{b}_{nv}\ \hat{a}'_{vr} \lambda_{v}
\right|^{k} \nonumber \\ & = & O(1) \sum_{n=2}^{\infty}n^{k-1} \left| \sum_{r=0}^{n-2}
\bar{\Delta}x_{r} \frac{b_{nn}\lambda_{n}}{a_{nn}}
\right|^{k}.
\end{eqnarray*}
Then as in $T_{n}(1)$, we have that
\begin{eqnarray*}
\overline{T_{n}(2)} & = &  \sum_{r=0}^{n-2} n^{1-\frac{1}{k}}
\bar{\Delta}x_{r} \frac{b_{nn}|\lambda_{n}|}{a_{nn}}    \nonumber \\    & = & \sum_{r=1}^{\infty} d_{nr} \bar{\Delta}x_{r}
\end{eqnarray*}
where
\[d_{nr}=\left\{\begin{array}{cl}
                                   n^{1-\frac{1}{k}} \frac{b_{nn}\lambda_{n}}{a_{nn}}  &, \mbox{    if $0\leq r\leq n-2$}\\
                                   0 &, \mbox{    if $r>n-2.$}
                                   \end{array}
        \right. \]
Now
\begin{eqnarray*}
\sum |\overline{T_{n}(2)}|^{k}<\infty \ \ whenever \ \ \sum |\bar{\Delta}x_{n}|<\infty
\end{eqnarray*}
is equivalently
\begin{eqnarray}
\sup_{r} \sum _{n=1}^{\infty} |d_{nr}|^{k}<\infty
\end{eqnarray}
by Lemma. But (26) is equivalent to
\begin{eqnarray}
\sum _{n=r}^{\infty} |d_{nr}|^{k} = O(1) \sum _{n=r+2}^{\infty} \left| n^{1-\frac{1}{k}} \frac{b_{nn}\lambda_{n}}{a_{nn}} \right|^{k} =O(1).
\end{eqnarray}
Therefore, we have
\begin{eqnarray*}
\sum_{n=1}^{\infty}n^{k-1}\left| T_{n}(i) \right| ^{k}<\infty
\quad for \quad i=1,2.
\end{eqnarray*}  
This completes the proof of theorem.

\end{document}